\documentstyle[11pt]{article}

\begin{document}

\def \qed {\hspace*{\fill}\frame{\rule[0pt]{0pt}{8pt}\rule[0pt]{8pt}{0pt}}\par}
\def \qedup{\vskip-20pt\qed}

\def\R{{\rm I\!R}} 
\def\N{{\rm I\!N}} 
\def\F{{\rm I\!F}}
\def\M{{\rm I\!M}}
\def\H{{\rm I\!H}}
\def\K{{\rm I\!K}}
\def\P{{\rm I\!P}}
\def\E{{\mathchoice {\rm 1\mskip-4mu l} {\rm 1\mskip-4mu l}
{\rm 1\mskip-4.5mu l} {\rm 1\mskip-5mu l}}}
\def\C{{\mathchoice {\setbox0=\hbox{$\displaystyle\rm C$}\hbox{\hbox
to0pt{\kern0.4\wd0\vrule height0.9\ht0\hss}\box0}}
{\setbox0=\hbox{$\textstyle\rm C$}\hbox{\hbox
to0pt{\kern0.4\wd0\vrule height0.9\ht0\hss}\box0}}
{\setbox0=\hbox{$\scriptstyle\rm C$}\hbox{\hbox
to0pt{\kern0.4\wd0\vrule height0.9\ht0\hss}\box0}}
{\setbox0=\hbox{$\scriptscriptstyle\rm C$}\hbox{\hbox
to0pt{\kern0.4\wd0\vrule height0.9\ht0\hss}\box0}}}}
\def\bbbe{{\mathchoice {\setbox0=\hbox{\smalletextfont e}\hbox{\raise
0.1\ht0\hbox to0pt{\kern0.4\wd0\vrule width0.3pt height0.7\ht0\hss}\box0}}
{\setbox0=\hbox{\smalletextfont e}\hbox{\raise
0.1\ht0\hbox to0pt{\kern0.4\wd0\vrule width0.3pt height0.7\ht0\hss}\box0}}
{\setbox0=\hbox{\smallescriptfont e}\hbox{\raise
0.1\ht0\hbox to0pt{\kern0.5\wd0\vrule width0.2pt height0.7\ht0\hss}\box0}}
{\setbox0=\hbox{\smallescriptscriptfont e}\hbox{\raise
0.1\ht0\hbox to0pt{\kern0.4\wd0\vrule width0.2pt height0.7\ht0\hss}\box0}}}}
\def\Q{{\mathchoice {\setbox0=\hbox{$\displaystyle\rm Q$}\hbox{\raise
0.15\ht0\hbox to0pt{\kern0.4\wd0\vrule height0.8\ht0\hss}\box0}}
{\setbox0=\hbox{$\textstyle\rm Q$}\hbox{\raise
0.15\ht0\hbox to0pt{\kern0.4\wd0\vrule height0.8\ht0\hss}\box0}}
{\setbox0=\hbox{$\scriptstyle\rm Q$}\hbox{\raise
0.15\ht0\hbox to0pt{\kern0.4\wd0\vrule height0.7\ht0\hss}\box0}}
{\setbox0=\hbox{$\scriptscriptstyle\rm Q$}\hbox{\raise
0.15\ht0\hbox to0pt{\kern0.4\wd0\vrule height0.7\ht0\hss}\box0}}}}
\def\T{{\mathchoice {\setbox0=\hbox{$\displaystyle\rm
T$}\hbox{\hbox to0pt{\kern0.3\wd0\vrule height0.9\ht0\hss}\box0}}
{\setbox0=\hbox{$\textstyle\rm T$}\hbox{\hbox
to0pt{\kern0.3\wd0\vrule height0.9\ht0\hss}\box0}}
{\setbox0=\hbox{$\scriptstyle\rm T$}\hbox{\hbox
to0pt{\kern0.3\wd0\vrule height0.9\ht0\hss}\box0}}
{\setbox0=\hbox{$\scriptscriptstyle\rm T$}\hbox{\hbox
to0pt{\kern0.3\wd0\vrule height0.9\ht0\hss}\box0}}}}
\def\bbqS{{\mathchoice
{\setbox0=\hbox{$\displaystyle     \rm S$}\hbox{\raise0.5\ht0\hbox
to0pt{\kern0.35\wd0\vrule height0.45\ht0\hss}\hbox
to0pt{\kern0.55\wd0\vrule height0.5\ht0\hss}\box0}}
{\setbox0=\hbox{$\textstyle        \rm S$}\hbox{\raise0.5\ht0\hbox
to0pt{\kern0.35\wd0\vrule height0.45\ht0\hss}\hbox
to0pt{\kern0.55\wd0\vrule height0.5\ht0\hss}\box0}}
{\setbox0=\hbox{$\scriptstyle      \rm S$}\hbox{\raise0.5\ht0\hbox
to0pt{\kern0.35\wd0\vrule height0.45\ht0\hss}\raise0.05\ht0\hbox
to0pt{\kern0.5\wd0\vrule height0.45\ht0\hss}\box0}}
{\setbox0=\hbox{$\scriptscriptstyle\rm S$}\hbox{\raise0.5\ht0\hbox
to0pt{\kern0.4\wd0\vrule height0.45\ht0\hss}\raise0.05\ht0\hbox
to0pt{\kern0.55\wd0\vrule height0.45\ht0\hss}\box0}}}}
\def\Z{{\mathchoice {\hbox{$\sf\textstyle Z\kern-0.4em Z$}}
{\hbox{$\sf\textstyle Z\kern-0.4em Z$}}
{\hbox{$\sf\scriptstyle Z\kern-0.3em Z$}}
{\hbox{$\sf\scriptscriptstyle Z\kern-0.2em Z$}}}}

\def\OO{{\cal O}}

\newcommand{\vv}{\vspace{0.3cm}} 
\newcommand{\vV}{\vspace{0.5cm}} 
\newcommand{\VV}{\vspace{2cm}} 

\title{Cusps and Codes}
\author{Wolf P. Barth \and S{\l}awomir Rams}
\date{{}}
\maketitle

\begin{abstract}
We study a construction, which produces surfaces
$Y \subset \P_3(\C)$ with cusps. For example we obtain
surfaces of degree six with 18, 24 or 27 three-divisible
cusps. For sextic surfaces in a particular family of
up to 30 cusps the codes of these sets of cusps
are determined explicitly.
\end{abstract}

\def\thefootnote{}
\footnotetext{
\hspace{-5.5ex} 
Suported by the DFG Schwerpunktprogramm
"Global methods in complex geometry". The second
author is supported by a
Fellowship of the Foundation for Polish Science
and KBN Grant No. 2 P03A 016 25. }

\newtheorem{defi}{Definition}[section]
\newtheorem{prop}{Proposition}[section]
\newtheorem{theo}{Theorem}[section]
\newtheorem{lemm}{Lemma}[section]

\vV

\setcounter{section}{-1}

\section{Introduction}

The aim of this note is to give examples of algebraic surfaces in
$\P_3(\C)$ with cusps, and to determine the code
of the set of these cusps.

Recall that a cusp (= singularity $A_2$) is a surface singularity
given in (local analytic) coordinates $x,y$ and $z$ 
centered at the singularity by an equation 
$$xy-z^3=0.$$
It is resolved by introducing two $(-2)$-curves. Let $Y \subset \P_3$
be an algebraic surface with $n$ cusps $P_1,...,P_n \in Y$. Let
$\pi:X \to Y$ be its minimal desingularisation with $E_{\nu}',E_{\nu}''$
the two $(-2)$-curves over $P_{\nu}$. The code of this set of cusps
is the kernel of the $\F_3$-linear morphism
$$\F_3^n \to H^2(X,\F_3), \quad 
(i_1,...,i_n) \mapsto \sum_{\nu=1}^n i_{\nu}[E_{\nu}'-E_{\nu}''].$$
A word $(i_1,...,i_n)$ belongs to this code if and only if the
class of the divisor $\sum i_{\nu} (E_{\nu}'-E_{\nu}'')$ is divisible
by 3 in $NS(X)$. Or equivalently: There is a cyclic triple
cover of $Y$ branched precisely over the points $P_{\nu}$ with
$i_{\nu} \not= 0$. Such a set of cusps is called 
3-divisible [B, T].

From
$$\left( \sum_{\nu=1}^n i_{\nu} [E_{\nu}'-E_{\nu}'']\right)^2 =
-6 \cdot \mbox{ number of } \nu \mbox{ with }i_{\nu} \not= 0 $$
it easily follows that the number $n$ of cusps in a 3-divisible
set is a multiple of $3$. 
For the maximal number of 3-divisible cusps on a surface $Y$ of
given degree $d$ there seems to be no better upper bound than the famous
bounds of Miaoka [M]
$$n \leq \frac{1}{4}d(d-1)^2$$
or Varchenko [V] for the maximal number of cusps, 3-divisible or not.
S.-L. Tan [T] shows that a surface of degree $d$ with $3 \leq d \leq 5$ 
can have only
$n$ three-divisible cusps with
$$ \begin{array}{c|c|c|c}
d & 3 & 4 & 5 \\ \hline
n & 3 & 6 & 12, 15, 18 \\
\end{array} $$
In the accompanying note [BR] we show that the minimal number of cusps
in a 3-divisible set on a sextic surface is 18. It is not known whether there
is a quintic surface with a 3-divisible set of 18 cusps, nor seems it to be 
known
for which numbers $n, \, 18 \leq n \leq 36,$ there is a sextic surface
with $n$ three-divisible cusps. Here we construct sextic surfaces with
$n=18,24$ and $27$ three-divisible cusps.

\vv
{\em Notations and conventions:} All varieties will be defined over the 
base-field $\C$. As coefficient field for cohomology we usually use
the field $\F_3$ with three elements. For brevity we denote it by $\F$.

\vV

\section{Constructions}

\subsection{The direct construction}

The basic idea for constructing surfaces $Y \subset \P_3$ with cusps
is very simple: just globalize the local equation $xy-z^3=0$. This means
the following: For fixed degree $d=6, \, 9, \, ...$, 
divisible by $3$, take polynomials
$s_1,...,s_k$ and $s$ of degrees
$$deg(s_1)=d_1,...,deg(s_k)=d_k, \quad d_1+...+d_k=d, \quad
deg(s)=d/3,$$
and consider the surface $Y \subset \P_3$ of degree $d$ given by the equation
$$s_1 \cdot ... \cdot s_k -s^3 = 0.$$
Let $S_1,...,S_k,S \subset \P_3$ be the surfaces with equation 
$s_1=0,...,s_k=0,s=0$. A Bertini-type-argument shows that for general choice
of these polynomials 

\begin{itemize}
\item any three surfaces $S_i,S_j,S$ meet transversally in 
$d_i \cdot d_j \cdot d/3$ points $P_{\nu}$, 
which then are cusps on the surface $Y$;
\item no four surfaces $S_i,S_j,S_m,S$ meet;
\item the surface $Y$ is smooth away from the cusps $P_{\nu}$ at 
the intersections
$S_i \cap S_j \cap S$.
\end{itemize}

The simplest examples of surfaces $Y$ with $n$ cusps obtained in this way are
for $d=6$
$$ \begin{array}{c|cccccccccc}
d_1,...,d_k & 1,5 & 2,4 & 3,3 & 1,1,4 & 1,2,3 & 2,2,2 & 1,1,1,3 & 1,1,2,2 & 1,1,1,1,2 & 1,1,1,1,1,1 \\
\hline
n         & 10 & 16 & 18 & 18 & 22  & 24  & 24   & 26   & 28    & 30     \\
\end{array}$$
and for $d=9$
$$ \begin{array}{c|cccccccccccc}
d_1,...,d_k & 1,8 & 2,7 & 1,1,7 & 3,6 & 4,5 & 1,2,6 & 1,1,1,6 & 1,3,5 & 1,4,4 & 
              2,2,5 & 1,1,2,5 & 2,3,4 \\ \hline  
n         & 24 & 42 & 45 & 54 & 60 & 60  & 63 & 69 &  72  & 72  & 75 & 78 \\ 
            \end{array} $$  
$$ \begin{array}{cccccccc} 
1,1,1,1,5 & 3,3,3 & 1,1,3,4 & 1,2,2,4 & 1,2,3,3 & 1,1,1,2,4 & 2,2,2,3 
                                                          & 1,1,1,3,3 \\ \hline
 78       & 81    &      81 &   84   & 87       &   87      & 90  & 90   \\
             \end{array} $$ 
$$\begin{array}{ccccc}
  1,1,1,1,1,4 & 1,1,1,1,2,3 & 1,1,1,1,1,1,3 &
  1,1,1,1,1,1,1,2 & 1,1,1,1,1,1,1,1,1 \\ \hline
  90     & 96    & 99 & 105 & 108  \\
\end{array}$$

\vV
\subsection{The residual construction}

The construction from section 1.1 produces only surfaces of degrees $d$ 
divisible by $3$. However in [BR] we observe:

\begin{itemize}
\item If $Y \subset \P_3$ is a quartic surface with 6 three-divisible cusps,
then there is a residual quadric $R$,
\item if $Y \subset \P_3$ is a quintic surface with 12 three-divisible
cusps, then there is a residual plane $R$,
\end{itemize}

\noindent
such that $Y \cup R$ has an equation $s_1 \cdot s_2-s^3=0$ with 
$deg(s_1)=deg(s_2)=3, \, deg(s)=2$. Here we use a residual 
surface $R$ to construct surfaces with a set of $3$-divisible cusps
in degrees $d$ not necessarily divisible by $3$. The Bertini-type
arguments for this 'residual construction' unfortunately are
more sophisticated than the argument needed in section 1.1. So we
restrict to the simplest case of the construction, which
probably can be generalized in several ways.
We start with

\begin{itemize}
\item a smooth residual surface $R: \, r=0$ of degree $b$;
\item $k$ auxiliary polynomials $r_1,...,r_k$ of degrees $c_1,...,c_k$
with $c_1+...+c_k=:c$
such that the curves $R_i: \, r=r_i=0$ on $R$ are smooth and intersect
transversally;
\item surfaces $S_i: \, s_i=0$ with $s_i=r_i^3+r \cdot t_i$ of degrees 
$d_i=3 c_i$;
\item a degree-$c$ surface $S: \, s=0$ with 
$s=r_1 \cdot ... \cdot r_k + r \cdot t.$
\end{itemize}

\noindent
Here we assume
$$ d_i \geq b, \quad c \geq b.$$
Then always
$$ s_1 \cdot ... \cdot s_k - s^3 = 
r \cdot (t_1 \cdot r_2^3 \cdot ... \cdot r_k^3 +...+
          r_1^3 \cdot ... \cdot r_{k-1}^3 \cdot t_k
         -3 \cdot r_1^2 \cdot ... \cdot r_k^2 \cdot t)\quad mod \quad r^2$$
vanishes on $R$. 
We are interested in the surface $Y$ defined by the polynomial
$$f:=\frac{1}{r} (s_1 \cdot ... \cdot s_k - s^3)$$
of degree $d=3c-b$ when the polynomials $t_i,t$ of degrees $d_i-b$
and $c-b$ are chosen generally.

1) Each polynomial $s_i$ vanishes along the curve $R_i$ to the first
order, except for points, where $t_i=0$. For general choice of $t_i$
this does not happen at an intersection $R_i \cap R_j$. And $s_i$ vanishes
to the second order at points on $R$, where $r_i=t_i=0$. So $f$ vanishes
on $R_i$ only at points with $r_i=r_j=0$ or $r_i=t_i=0$, and there
to the first order. This shows that $Y$ is smooth wherever it meets
any curve $R_i$.

2) The polynomials 
$$f|R = (t_1 \cdot r_2^3 \cdot ... \cdot r_k^3+ ...
+r_1^3 \cdot ... \cdot r_{k-1}^3 \cdot t_k 
-3 \cdot r_1^2 \cdot ... \cdot r_k^2 \cdot t)|R, \quad
t_1,...,t_k,t \mbox{ varying}$$
form a linear system with base locus consisting of the finitely many
points $r=r_i=r_j=0$. At these points $df \not= 0$ by 1). So for general
choice of $t_1,...,t_k,t$ the surface $f=0$ is smooth at its 
intersection with $R$.

3) The polynomials $s = r_1 \cdot ... \cdot r_k + r \cdot t$ 
with $t$ vaying form a linear 
system with base locus $R_1 \cup ... \cup R_k$. 
If $t$ does not vanish at any intersection $R_i \cap R_j$ the there
$ds \not= 0$.
So for general choice of $t$ the 
surface $S$ is smooth.

4) For each $i$ the polynomials $s_i=r_i^3+r \cdot t_i$ form a linear system
with base locus $R_i$. So for general choice of $t_i$ the surface $S_i: s_i=0$
as well as the curve $C_i: \, s_i=s=0$
is smooth outside of $R$.
And for general choice of $t_i$ and $t_j$ the
surfaces $S_i,S_j,S$ intersect transversally outside of $R$.

5) The polynomials 
$$(r_1^3+r \cdot t_1)\cdot s_2 \cdot ... \cdot s_k-s^3 =
r_1^3 \cdot s_2 \cdot ... \cdot s_k-s^3 +t_1 \cdot 
r \cdot s_2\cdot ... \cdot s_k, \quad t_1 \mbox{ varying}$$
form a linear system with base locus $Y \cap R$, where $Y$ is smooth by 2),
and the curves $s=s_2=0,...,s=s_k=0$. So for general choice of $t_1$ the
surface $s_1 \cdot ... \cdot s_k-s^3$ is smooth outside of $S$.

6) On $S_i \cap S$ we have
$$d(s_1 \cdot ... \cdot s_k-s^3) = 
\sum_i s_i \cdot \prod_{j \not= i} s_j.$$
By 4) this differential is $\not= 0$ outside of $R$ and away
from the points in $S_i \cap S_j, \, j \not= i$.

Statements 1)-6) show that the surface $Y: \, f=0$ is smooth away from
the points in $S_i \cap S_j \cap S=0$. By 4) the three surfaces 
$S_i,S_j,S$ intersect transversally in these points, hence 
these points are cusps 
on $Y$.

\begin{lemm} In each point $P \in S_i \cap S_j \cap R$ we have the  
intersection number  
$$i_P(S_i,S_j,S)=6.$$ \end{lemm}

Proof. By 4) we may take $x=r_i,y=r_j,z=r$ as local coordinates.
Since $S: \, xyg(x,y,z)+zt=0$ with $g(0,0,0) \not= 0$ 
is smooth, we may eliminate $z=-xyg/t$. The 
intersection number then is the intersection number of the two curves
$$s_i=x \cdot (x^2-\frac{t_ig}{t}y)=0, \quad
s_j=y \cdot (y^2-\frac{t_jg}{t}x)=0.$$ \qedup     

Altogether we found:

\begin{theo} The residual construction with the degrees taken as above gives 
surfaces $Y$ of degree $d=3c - b$ which are smooth but for 
$$n_{i,j} = 3c_i \cdot 3c_j \cdot c - 6 \cdot c_i \cdot c_j \cdot b
=3 \cdot c_i \cdot c_j \cdot (d - b)$$
cusps on each curve $S_i \cap S_j$. \end{theo} 

In the following table we give some examples.

$$ \begin{array}{rrrrr|c|ccc}
d_1 & d_2 & d_3 & c & b & d & n_{1,2} & n_{1,3} & n_{2,3} \\ \hline
3   & 3   &  -  & 2 & 2 & 4 &    6    &  -      &    -    \\ 
3   & 3   &  -  & 2 & 1 & 5 &    12   &  -      &    -    \\ \hline
3   & 6   &  -  & 3 & 3 & 6 &    18   &  -      &    -    \\
3   & 6   &  -  & 3 & 2 & 7 &    30   &  -      &    -    \\
3   & 6   &  -  & 3 & 1 & 8 &    42   &  -      &    -    \\ \hline
3   & 3   &  3  & 3 & 3 & 6 &     9   &  9      &    9    \\
3   & 3   &  3  & 3 & 2 & 7 &    15   &  15     &   15    \\
3   & 3   &  3  & 3 & 1 & 8 &    21   &  21     &   21    \\
\end{array} $$

\vV

\subsection{Codes}

Denote by $C_i =R_i \cup B_i$ the curve $s_i=s=0$. For general choice
of the polynomials above, $B_i$ is smooth away from $R$. 
We have $S_i \cap Y = B_i$, both the surfaces $S_i$ and $Y$ touching
along $B_i$ to the third order.
Let $\pi:X \to Y$ be the minimal resolution of $Y$, obtained e.g.
by blowing up $\P_3$ in the cusps $P_{\nu}$ of $Y$. Over each cusp 
$P_{\nu}$ it introduces two $(-2)$-curves $E_{\nu}'$ and $E_{\nu}''$.
If $P_{\nu} \in B_i \cap B_j$, then the proper transform $D_i$ of
$B_i$ in $X$ meets (transversally) exactly one of the two exceptional
curves over $P_{\nu}$, while the proper transform $D_j$ of $B_j$ meets the 
other
exceptional curve.

By abuse of notation we put 
$$\OO_X(m):= \pi^* \OO_Y(m), \quad m \in \Z.$$
Let $P_{\nu}$ be the cusps of $Y$ in $B_i$. Let 
$E_{\nu}' \subset X$ be the $(-2)$-curve over $P_{\nu}$ meeting $D_i$
while $E_{\nu}''$ meets some curve $D_j, \, j \not= i$. Then
$$\OO_X(d_i) \sim 3D_i + \sum_{\nu} (i_{\nu}' E_{\nu}'+i_{\nu}''E_{\nu}'')$$
with integers $i_{\nu}',i_{\nu}'' \geq 0$. Here 
$(\OO_X(1).E_{\mu}')=(\OO_X(1).E_{\mu}'')=0$ implies
$$ (3D_i+\sum (i_{\nu}'E_{\nu}'+i_{\nu}'' E_{\nu}'')).E_{\mu}'
=3 -2i_{\mu}'+i_{\mu}'' = 0, $$
$$ (3D_i+\sum (i_{\nu}'E_{\nu}'+i_{\nu}'' E_{\nu}'')).E_{\mu}''
=i_{\mu}'-2i_{\mu}'' = 0. $$
For $P_{\nu} \in B_i$ we find
$$i_{\nu}'=2i_{\nu}'', \quad i_{\nu}''=1, i_{\nu}'=2,$$
and the class
$$\OO_X(d_i) - \sum_{P_{\nu} \in C_i} (2E_{\nu}'+E_{\nu}'')
\sim 3 D_i \in NS(X)$$
is divisible by $3$.

\begin{defi} The {\em code} of $Y$ is the kernel of the $\F$-linear
morphism
$$\F^n \to H^2(X , \F), \quad
(i_1,...,i_n) \mapsto \sum_{\nu} i_{\nu}[E_{\nu}'-E_{\nu}''].$$
The {\em extended code [E]} is the kernel of
$$\F^{n+1} \to H^2(X , \F), \quad
(i_0, i_1,...,i_n) \mapsto 
\OO_X(i_0)+ \sum_{\nu} i_{\nu}[E_{\nu}'-E_{\nu}''].$$
\end{defi}

To avoid a clumsy notation we put
$$ e_{\nu} := [E_{\nu}'-E_{\nu}''] = -[2E_{\nu}'+E_{\nu}''] \,\,
mod \, \, 3,$$
where the exceptional curves are ordered such that
$$3 D_i \sim \OO_X(d_i) + 
\sum_{P_{\nu} \in C_i \cap C_j, \, j<1} -e_{\nu} 
+\sum_{P_{\nu} \in C_i \cap C_j, \, j>i} e_{\nu}.$$
For each $i$ we then obtain a word
$$w_i=(d_i \, mod \, 3, \, \underbrace{-1,...,-1}_{j<i}, \, 
               \underbrace{1,...,1}_{j>i} )$$
in the extended code. 

\VV
\section{Sextics}

Here we determine the codes of the sextic surfaces $Y$ given
by the construction in section 1.1. 
We fix a partition $d_1,...,d_k$ of $6$. An 
equation 
$$f:=s_1 \cdot ... \cdot s_k-s^3 = 0, \quad deg(s_i)=d_i, \, deg(s)=2$$
will be called an equation of {\em type} $d_1,...,d_k$. 
It will be called {\em admissible}, if

\begin{itemize}
\item Any three surfaces $S_i,S_j,S$ meet transversally at
$2 \cdot d_i \cdot d_j$ points. These points then are cusps on $Y$.
\item There are no other singularities on $Y$ but these cusps.
\end{itemize}

A Bertini-Argument shows that for each partition of $6$ 
admissible equations of this type exist. 
Then they are dense in the family of all equations of a given type.
All equations of a fixed type form an irreducible family. Hence
all admissible equations of given type form a connected family.
Each path $f_t, \, 0 \leq t \leq 1,$ with $Y_0: \, f_0=0, \,
Y_1: \, f_1=0$ defines an isomorphism 
$H^2(X_0,\F) \to H^2(X_1,\F)$ inducing an isomorphism of
codes. This shows

\begin{prop} All sextic surfaces with admissible 
equation of a fixed type have isomorphic codes.
\end{prop}

\begin{defi} The proper code $C_{d_1...d_k}$, resp. the extended code 
$E_{d_1...d_k}$ of type $d_1,...,d_k$ 
is the proper, resp. extended code of all the surfaces having an 
admissible equation of type $d_1,...,d_k$. \end{defi}

\vv
To understand the codes of our sextic surfaces we use

\begin{prop} {\bf a)} Each word in the proper code has a weight $n\geq 18$.

\noindent
{\bf b)} If the partition $c_1,...,c_l$ of $6$ is a sub-partition of
$d_1,...,d_k$, then the code of type $d_1,...,d_k$ is isomorphic
to a subcode of the code of type $c_1,...,c_l$.

\noindent
{\bf c)} All codes of our surfaces admit an involution, which arises from
interchanging the cusps in pairs. \end{prop}

Proof. a) The lower bound $n \geq 18$ for $n$ three-divisible cusps on a sextic
surface is proven in [BR, thm.1.1]. 

b)  It suffices to prove the assertion for $c_1,...,c_{k+1}=
d_1,...,d_{k-1},c_k,c_{k+1}$ with $d_k=c_k+c_{k+1}$. 
So fix some admissible equation $s_1 \cdot ... \cdot  s_k = s^3$ of type
$d_1,...,d_k$ and an admissible equation 
$s_1 \cdot ... \cdot s_{k-1} \cdot t_k \cdot t_{k+1} = s^3$ of type
$d_1,...,d_{k-1},c_k,c_{k+1}$. Consider the one-parameter-family
of surfaces with equation
$$s_1 \cdot ... \cdot s_{k-1} \cdot (\lambda s_k + 
(1- \lambda)\cdot t_k \cdot t_{k+1}) = s^3, \quad \lambda \in \C.$$
Define the surfaces $T_k: \, t_k=0$ and $T_{k+1}: \, t_{k+1}=0$.
For all $\lambda \not= 0$, but finitely many, the equation is
admissible of type $d_1,...,d_k$. For all $i,j \leq k-1$
the cusps $P_{\mu} \in S_i \cap S_j \cap S$ coincide. And for
$i \leq k-1$ the cusps $P_{\nu}$ on $S_i \cap S_k \cap S$
converge for $\lambda \to 0$ to cusps on $S_i \cap T_k \cap S$
or $S_i \cap T_{k+1} \cap S$. This shows that the code
of type $d_1,...,d_k$ is the subcode of the code of type
$d_1,...,d_{k-1},c_k,c_{k+1}$ which consists of the words
assigning the value $0$ to the cusps in $T_k \cap T_{k+1} \cap S$.

c)  Fix the quadric $S: \, s:=x_2x_3=0$ 
admitting the involution 
$$I: (x_0:x_1:x_2:x_3) \mapsto (x_0:x_1:x_3:x_2),$$
which interchanges the two planes of $S$.
Choose $f_1(x_0,x_1,x_2),...,f_k(x_0,x_1,x_2)$ of degrees
$i_1,...,i_k$ such that all the curves $C_i: f_i=0$ on 
the plane $x_3=0$ are smooth and intersect
transversally, not on the line $x_2=x_3=0$. 
Then put $s_i(x_0,...,x_3):=f_i(x_0,x_1,x_2+x_3), \, i=1,...,k$. 
The polynomials $s_i$ are $I$-invariant and define surfaces
$S_i:\, s_i=0$ meeting transversally on $S$. The surface
$s_1 \cdot ... \cdot s_k-s^3$ has cusps in the points on
$S$ where $s_i = s_j=0$. Elsewhere it is smooth along its
intersection with $S$. So by Bertini, 
for general $\lambda$
the $I$-invariant 
surface of equation $\lambda \cdot s_1 \cdot ... \cdot s_k = s^3$
is smooth but for these points on $S$. \qed

\vv
In the following table we give words, which by 1.3 
belong to the extended code. A pair $ij$ in the first row stands for
the cusps in the intersection $S_i \cap S_j$. For each type
the number below $ij$ is the number of cusps in the corresponding
intersection. The words are specified by their values taken at these cusps.

$$ \begin{array}{l|l|rrrrrrrrrrrrrrr}
\mbox{type} &i_0&12&13&23&14&24&34&15&25&35&45&16&26&36&46&56 \\  \hline
1,5         &   &10&  &  &  &  &  &  &  &  &  &  &  &  &  &   \\  
w           &1  & 1&  &  &  &  &  &  &  &  &  &  &  &  &  &   \\ \hline 
2,4         &   &16&  &  &  &  &  &  &  &  &  &  &  &  &  &   \\  
w           &2  & 1&  &  &  &  &  &  &  &  &  &  &  &  &  &   \\ \hline 
3,3         &   &18&  &  &  &  &  &  &  &  &  &  &  &  &  &   \\  
w           &0  & 1&  &  &  &  &  &  &  &  &  &  &  &  &  &   \\ \hline 
1,1,4       &   & 2& 8& 8&  &  &  &  &  &  &  &  &  &  &  &   \\ \hline 
w_1         &1  & 1& 1&  &  &  &  &  &  &  &  &  &  &  &  &   \\  
w_2         &1  & 2&  & 1&  &  &  &  &  &  &  &  &  &  &  &   \\ \hline 
1,2,3       &   & 4& 6&12&  &  &  &  &  &  &  &  &  &  &  &   \\ \hline 
w_1         &1  & 1& 1&  &  &  &  &  &  &  &  &  &  &  &  &   \\  
w_2         &1  & 2&  & 1&  &  &  &  &  &  &  &  &  &  &  &   \\ \hline 
2,2,2       &   & 8& 8& 8&  &  &  &  &  &  &  &  &  &  &  &   \\ \hline 
w_1         &1  & 1& 1&  &  &  &  &  &  &  &  &  &  &  &  &   \\  
w_2         &1  & 2&  & 1&  &  &  &  &  &  &  &  &  &  &  &   \\ \hline 
1,1,1,3     &   & 2& 2& 2& 6& 6& 6&  &  &  &  &  &  &  &  &   \\ \hline 
w_1         &1  & 1& 1&  & 1&  &  &  &  &  &  &  &  &  &  &   \\  
w_2         &1  & 2&  & 1&  & 1&  &  &  &  &  &  &  &  &  &   \\ 
w_3         &1  &  & 2& 2&  &  & 1&  &  &  &  &  &  &  &  &   \\ \hline 
1,1,2,2     &   & 2& 4& 4& 4& 4& 8&  &  &  &  &  &  &  &  &   \\ \hline 
w_1         &1  & 1& 1&  & 1&  &  &  &  &  &  &  &  &  &  &   \\  
w_2         &1  & 2&  & 1&  & 1&  &  &  &  &  &  &  &  &  &   \\ 
w_3         &2  &  & 2& 2&  &  & 1&  &  &  &  &  &  &  &  &   \\ \hline 
1,1,1,1,2   &   & 2& 2& 2& 2& 2& 2& 4& 4& 4& 4&  &  &  &  &  \\ \hline 
w_1         &1  & 1& 1&  & 1&  &  & 1&  &  &  &  &  &  &  &   \\  
w_2         &1  & 2&  & 1&  & 1&  &  & 1&  &  &  &  &  &  &   \\ 
w_3         &1  &  & 2& 2&  &  & 1&  &  & 1&  &  &  &  &  &   \\  
w_4         &1  &  &  &  & 2& 2& 2&  &  &  & 1&  &  &  &  &   \\ \hline 
1,1,1,1,1,1 &   & 2& 2& 2& 2& 2& 2& 2& 2& 2& 2& 2& 2& 2& 2& 2 \\ \hline 
w_1         &1  & 1& 1&  & 1&  &  & 1&  &  &  & 1&  &  &  &   \\  
w_2         &1  & 2&  & 1&  & 1&  &  & 1&  &  &  & 1&  &  &   \\ 
w_3         &1  &  & 2& 2&  &  & 1&  &  & 1&  &  &  & 1&  &   \\  
w_4         &1  &  &  &  & 2& 2& 2&  &  &  & 1&  &  &  & 1&   \\  
w_5         &1  &  &  &  &  &  & 2& 2& 2& 2&  &  &  &  &  & 1 \\  
\end{array} $$

All words given in this table are linearly independent in $E_{d_1...d_k}$.

\vV
\subsection{The case 1,1,1,1,1,1}

We fix an involution $I$ as in prop. 2.2 c) and order the cusps
$P_{\nu}$ such that $I$ interchanges $P_{\nu}$ and $P_{15+\nu}$ for 
$\nu=1,...,15$. This induces an involution on $\{1,...,30\}$ and on
$\F^{30}$. The code  $C_{1,1,1,1,1,1} \subset \F^{30}$ is invariant
under this involution. So it splits as a direct sum $C^+ \oplus C^-$ with
$C^+$ consisting of invariant words and $C^-$ of anti-invariant ones.

The symmetric group $\Sigma_6$ 
acts on the set of six
planes $S_i$ by permutations. 
Such a permutation can be realized by a path in the space of all admissible
equations. This shows that there is a $\Sigma_6$-action on the code.

\begin{prop} The actions of $I$ and $\Sigma_6$ can be chosen such that they 
commute. \end{prop}

Proof. As in prop. 2.2 consider the quadric $S: x_2 \cdot x_3=0$  
and the involution $I:x_2 \leftrightarrow x_3$.
Let $P_1,...,P_{15}$ be the cusps in the plane $x_2=0$ and let 
$P_{16},...,P_{30}$
be their images in the plane $x_3=0$. So $I(P_{\nu})=P_{15+\nu}$ for
$\nu=1,...,15$. 
Following a path in the space of
admissible equations does not interchange cusps between the two planes.
So $I \sigma (\nu)=\sigma I (\nu)$ for all $\nu$ and $\sigma \in \Sigma_6$. 
This is
the assertion. \qed

\vv
The cusps $P_1,...,P_{15}$ can be relabelled $P_{i,j}, \, i<j \leq 6$
such that $P_{i,j} \in S_i \cap S_j$. This induces an identification
$\F^{15} = \Lambda^2 (\F^6) \subset \F^{30}$. 
Let $e_1,...,e_6$ be canonical
generators for $\F^6$.
Each word $w_i \in E_{1,1,1,1,1,1}$ induces the word
$$u_i:=\sum_{k=1}^6 e_i \wedge e_k \in E_{111111} \cap \Lambda^2(\F^6).$$
For $i < j$ the word
$$u_{i,j} :=u_i - u_j = (e_i-e_j) \wedge (e_1+...+e_6)=
2e_i \wedge e_j + (e_i-e_j) \wedge \sum_{k \not= i,j} e_k$$
belongs to $C^+ \cap \Lambda^2(\F^6)$. 

Let $e=e_1+...+e_6 \in \F^6$ and $U \subset \F^6$
be the hyperplane of $u =\sum_i u_i e_i$ with $\sum u_i=0$. Notice
$e \in U$. 

\begin{lemm} The words $u_{i,j} \in \Lambda^2(\F^6)$ generate
the subspace $e \wedge U \subset \Lambda^2(\F^6)$. \end{lemm}

Proof. Obviously the words $u_{i,j}$ belong to $e \wedge U$.
This subspace has dimension $4$. And it is easy to see that the
four words $u_{1,2}, ..., u_{1,5}$ are linearly independent. \qed

\vv
We now put
$$B^+:=C^+ \cap \Lambda^2(\F^6), \quad
B^-:=C^- \cap \Lambda^2(\F^6).$$
Our main technical result is

\begin{prop} a) $B^+ \subset \Lambda^2(\F^6)$ coincides with $e \wedge U$. 
\hspace{1cm} b) $B^-\subset B^+$. \end{prop}

This implies our main result:

\begin{theo} The code $E_{111111}$ is generated by the words
$w_1,...,w_5$ from the table at the beginning of section 2. \end{theo}

Proof. Recall $C = C^+ \oplus C^-$. Now $C^+=\{(w,w), w \in B^+\}$.
By prop. 2.4 a)  $B^+$ is generated by the words $u_{i,j}$
and therefore $C^+$ is generated by the $w_i-w_j$. And
$C^-=\{(w,-w), \, w \in B^-\}$. Since $B^- \subset B^+$ by prop. 2.4 b),
for each word $(w,-w) \in C^-$ there is a word $(w,w) \in C^+$.
The word $(w,w)-(w,-w) \in C_{111111}$ has weight $\leq 15 < 18$.
So prop. 2.2 a)  implies $w=0$. We found $C^-=0$ and $C_{111111}=C^+$.
Then $E_{111111}$ is generated by the words $w_i$. \qed

\vv
Each code $B^+,B^- \subset \Lambda^2(\F^6)$ has the following two
properties:

\begin{itemize}
\item[1)] All its words have length $0,9,12,$ or $15$.
\item[2)] It is invariant under the action of $\Sigma_6$
induced from the permutations of coordinates in $\F^6$.
\end{itemize}

Proposition 2.4 therefore follows from

\begin{prop} Let $V \subset \Lambda^2(\F^6)$ be some code
with the properties above. Then $V \subset e \wedge U$. \end{prop}

Proof. We abbreviate $e_{ij}=e_i \wedge e_j =-e_{ji}$ for $i<j$.
We proceed in several steps.

{\em Step 1: Let $v=\sum_{i < j} v_{ij} e_{ij} \in V$ with
some $v_{ij}=0$. Then $v_{ik}=v_{jk}$ for all $k=1,...,6$.}

Proof. Write 
$$v =\sum_{k \not=i,j} (v_{ik}e_{ik}+v_{jk}e_{jk})
+\sum_{k,l \not= i,j} v_{kl}e_{kl}.$$
The involution $(i,j) \in \Sigma_6$ maps $v$ onto
$$(i,j)v =\sum_{k \not=i,j} (v_{jk}e_{ik}+v_{ik}e_{jk})
+\sum_{k,l \not= i,j} v_{kl}e_{kl}.$$
So
$$v-(i,j)v =\sum_{k \not=i,j} (v_{ik}-v_{jk})(e_{ik}-e_{jk})$$
has length $\leq 8$. By property $1)$ this word is $0$ and
$v_{ik}=v_{jk}$ for all $k$. \qed

{\em Step 2:} Given $v = \sum v_{ij} e_{ij} \in V$ consider
$u=v_{12}u_{12}+v_{34}u_{34}+v_{56}u_{56} \in e \wedge U$.  
Then $v'=v+u$ has the coefficients $v'_{12}=v'_{34}=v'_{56}=0$.
By step 1
$$v' = a \cdot (e_1+e_2) \wedge (e_3+e_4) 
      +b \cdot (e_1+e_2) \wedge (e_5+e_6)
      +c \cdot (e_3+e_4) \wedge (e_5+e_6), \quad a,b,c \in \F.$$
The assertion follows, if we show $v' \in e \wedge U$.

{\em Step 3:} We simplify the notation putting $v'=v$ and have
to show $v \in e \wedge U$. We apply symmetries from $\Sigma_6$ 
to obtain
\begin{eqnarray*}
(13)(24)v &=& -a \cdot (e_1+e_2) \wedge (e_3+e_4)
              +c \cdot (e_1+e_2) \wedge (e_4+e_5)
              +b \cdot (e_3+e_4) \wedge (e_4+e_5), \\
(15)(26)v &=& -c \cdot (e_1+e_2) \wedge (e_3+e_4)
              -b \cdot (e_1+e_2) \wedge (e_4+e_5)
              -a \cdot (e_3+e_4) \wedge (e_4+e_5).
\end{eqnarray*}
From $weight(v+(13)(24)v) \leq 8$ we conclude that $b=-c$  and
from $weight(v+(15)(26)v) \leq 8$ that $a=c$. The word
\begin{eqnarray*}
v &=& a \cdot ((e_1+e_2) \wedge (e_3+e_4)
              -(e_3+e_4) \wedge (e_5+e_6)
              +(e_3+e_4) \wedge (e_5+e_6)) \\
  &=& a \cdot (e_3+e_4-e_1-e_2) \wedge e
\end{eqnarray*}
therefore belongs to $e \wedge U$. \qed

\vV
\subsection{The other cases}

By prop. 2.2 b) 
there are inclusions of extended codes. The following table shows these 
inclusions, the bottom row giving the dimensions of the codes generated 
by the words $w_i$ from the table at the beginning of section 2:
$$ \begin{array}{ccccccccc}
33 &      &     &      &      &     &       &     &        \\
   & \searrow & &      &      &     &       &     &        \\
15 & \to  & 123 & \to  & 1113 & \to & 11112 & \to & 111111 \\
   & \searrow \hspace{-0.4cm} \nearrow & &
     \searrow \hspace{-0.4cm} \nearrow &  & \nearrow &  &     &        \\
24 & \to  & 114 & \to  & 1122 &     &       &     &        \\
   & \searrow & & \nearrow  & &     &       &     &        \\
   &      & 222 &      &      &     &       &     &        \\
   &      &     &      &      &     &       &     &        \\
1  &      & 2   &      &  3   &     &   4   &     &   5    \\
\end{array} $$
With each inclusion in this table new cusps appear, and the bigger code
contains a word taking non-zero values at the new cusps. So the bigger
code has a larger dimension than the included code. But the
right-hand code $E_{111111}$ has dimension $5$ by theorem 2.1. This proves

\begin{theo} For all sextic surfaces 
constructed by the method 1.1 the extended code
is generated by the words given in the table at the beginning of section 2. 
Its dimension is the
number in the table above. \end{theo} 

As a consequence we have:

\begin{theo} The two types $3,3$ and $1,1,4$ of surfaces with 
$18$ cusps, as well as the two types $2,2,2$ and $1,1,1,3$ of
surfaces with $24$ cusps differ by the dimensions of their
extended codes. So they cannot belong to the same connected family
of surfaces with their number of cusps, nor can they be 
degenerations of each other. \end{theo}

\vV
\subsection{27 cusps}

We apply the residual construction of section 1.2 for 
$c_1=c_2=c_3=1, \, b=3$ to obtain a sextic surface $Y$ with nine
cusps on each curve $S_i \cap S_j$. By section 1.3 the proper code contains
two words
$$ \begin{array}{c|ccc}
    & i_1,...,i_9 & i_{10},...,i_{18} & i_{19},...,i_{27} \\ \hline
w_1 &      0      &        1         &        1          \\
w_2 &      2      &        0         &        1          \\
\end{array} $$
and the word $w_1+w_2$ of weight $27$.

We finish with a quite explicit example. Take as residual cubic the Fermat
cubic defined by 
$$r:=x_0^3+x_1^3+x_2^3+x_3^3 $$
and put
$$s_i:=x_i^3+\lambda_i \cdot r, \, i=1,2,3, \quad s=x_1 \cdot x_2 \cdot x_3.$$
This leads to the sextic defined by
\begin{eqnarray*}
f &=& \frac{1}{r} ((x_1^3+\lambda_1r)(x_2^3+\lambda_2r)(x_3^3+\lambda_3r)
                    -x_1^2 \cdot x_2^3 \cdot x_3^3) \\
  &=& \lambda_1 \cdot x_2^3 x_3^3+\lambda_2 \cdot x_1^3 x_3^3
     +\lambda_3 \cdot x_1^3 x_2^3 +
      (\lambda_1 \lambda_2 \cdot x_3^3+\lambda_1 \lambda_3 \cdot x_2^3
     + \lambda_2 \lambda_3 \cdot x_1^3)\cdot r 
     + \lambda_1 \lambda_2 \lambda_3 \cdot r^2.
\end{eqnarray*} Here the surface $S: \, s=0$ of course is not smooth,
however it can be checked by direct computation that for general
choice of the coefficients $\lambda_i$ the sextic surface is smooth,
but for 27-cusps on the coordinate planes, like the nine points
in the intersection $x_1=s_2=s_3=0$. 

As far as we know, this number 27 is the largest number of a set of
3-divisible cusps on a sextic surface so far observed.

\VV

\noindent
Wolf P. Barth, 
Mathematisches Institut der Universit\"{a}t, Bismarckstr. 1 1/2,
D - 91054 Erlangen,
e-mail:barth@mi.uni-erlangen.de

\vspace{1ex}
\noindent
S{\l}awomir Rams, 
Institute of Mathematics, Jagiellonian University,
ul.~Reymonta~4, 30-059 Krak\'ow,
e-mail:rams@im.uj.edu.pl

\end{document}